\documentclass[11pt,a4paper]{amsart}

\title[One parameter regularizations of products of residue
currents]{One parameter regularizations of  products \\ of residue currents}

\author[M.\ Andersson \& H.\ Samuelsson Kalm \& E.\ Wulcan \& A.\
Yger]{Mats Andersson \& H{\aa}kan Samuelsson Kalm \&\\  Elizabeth Wulcan \& Alain Yger}
\thanks{First three authors partly supported by the Swedish Research
  Council.}
\address{M. Andersson, H. Samuelsson Kalm, E. Wulcan, Department of Mathematical Sciences, Division of Mathematics, 
University of Gothenburg and 
Chalmers University of Technology, SE-412 96 G\"{o}teborg, Sweden}
\email{matsa@chalmers.se, hasam@chalmers.se, wulcan@chalmers.se}
\address{A. Yger, Institut de Math\'{e}matiques, Universit\'{e} Bordeaux 1, 33405 Talence, France}
\email{Alain.Yger@math.u-bordeaux1.fr}
\date{\today}

\usepackage[english]{babel}
\usepackage{amsmath}
\usepackage{amsthm,amssymb,latexsym}
\usepackage{url,ae}
\usepackage{t1enc}
\usepackage{mathrsfs}
\usepackage{graphicx}
\usepackage{eufrak} 
\usepackage{geometry}
\newtheorem{proposition}{Proposition}[section]
\newtheorem{theorem}[proposition]{Theorem}
\newtheorem{lemma}[proposition]{Lemma}

\theoremstyle{definition}

\newtheorem{example}[proposition]{Example}
\newtheorem{remark}[proposition]{Remark}

\numberwithin{equation}{section}

\newcommand{\C}{\mathbb{C}}
\newcommand{\dbar}{\bar{\partial}}

\newcommand{\hol}{\mathscr{O}}
\newcommand{\w}{\wedge}

\def\newop#1{\expandafter\def\csname #1\endcsname{\mathop{\rm #1}\nolimits}}
\newop{span}

\begin{document}
\nocite{*}
\bibliographystyle{plain}

\begin{abstract}
We show that Coleff-Herrera type products of residue currents can be 
defined by analytic continuation of natural functions depending on 
one complex variable.
\end{abstract}

\maketitle
\thispagestyle{empty}

\section{Introduction}
Let $f$ be a holomorphic function defined on a domain in $\C^n$. 
It is proved in \cite{HL} using Hironaka's desingularization theorem that if $\varphi$ is a test form then  
\begin{equation*}
\lim_{\epsilon \to 0^+}\int_{|f|^2>\epsilon}\varphi/f
\end{equation*}
exists and defines the action of a current, denoted $1/f$. The $\dbar$-image, $\dbar (1/f)$, is the residue 
current of $f$ and it has the useful property that it is annihilated by a holomorphic function $g$ if and
only if $g$ is in the ideal generated by $f$.
If $f_1,\ldots,f_q$ are holomorphic functions then
the \emph{Coleff-Herrera product} of the currents $\dbar (1/f_j)$ is defined as follows. For a test form $\varphi$
of bidegree $(n,n-q)$ consider the residue integral
\begin{equation*}
I_f^{\varphi}(\epsilon) = \int_{T(\epsilon)} \frac{\varphi}{f_1\cdots f_q},
\end{equation*}
where $T(\epsilon)=\cap_1^q \{|f_j|^2=\epsilon_j\}$. It is proved in \cite{CH} that the limit of 
$\epsilon \mapsto I_f^{\varphi}(\epsilon)$ exists if $\epsilon=(\epsilon_1,\ldots,\epsilon_q)\to 0$ 
along a path in $\mathbb{R}_+^q$ such that $\epsilon_j/\epsilon_{j+1}^k\to 0$ for all $k\in \mathbb{N}$ and $j=1,\ldots,q-1$;
such a path is said to be \emph{admissible}. Moreover, the limit defines the action of a current, the Coleff-Herrera product
\begin{equation}\label{eq:CHprod}
\dbar \frac{1}{f_q}\wedge \cdots \wedge \dbar \frac{1}{f_1} . \, \varphi :=
``\lim_{\epsilon\to 0}" I_f^{\varphi}(\epsilon),
\end{equation}
where $``\lim"$ means the limit along an admissible path as above.
Following Passare \cite{PCrelle}, let $\chi$ be a smooth approximation of the characteristic function 
$\mathbf{1}_{[1,\infty)}$ and consider the smooth form
\begin{equation}\label{eq:passare}
\frac{\dbar \chi(|f_q|^2/\epsilon_q)}{f_q}\wedge \cdots \wedge \frac{\dbar \chi(|f_1|^2/\epsilon_1)}{f_1}.
\end{equation}
It follows from \cite[Theorem~2]{LS} or the proof of \cite[Proposition~2]{PCrelle} that the limit in the sense of currents
of \eqref{eq:passare} as $\epsilon \to 0$ along an admissible path equals the Coleff-Herrera product, and moreover,
that one gets the same result if one first lets $\epsilon_1\to 0$, then lets $\epsilon_2\to 0$ and so on. 
The Coleff-Herrera product 
is thus indeed the result of an iterative procedure. In general there are no obvious commutation properties, e.g.,
$\dbar (1/zw)\wedge \dbar (1/z)=0$ whereas $\dbar (1/z)\wedge \dbar (1/zw)=\dbar(1/z^2)\wedge \dbar (1/w)$, where
the last product is simply a tensor product.
However, if $f=(f_1,\ldots,f_q)$ defines a complete intersection, i.e., $\textrm{codim}\, \{f=0\} = q$,
then the Coleff-Herrera product depends in an anticommutative way of the ordering
of the tuple $f$; in fact by \cite{JEBHS} the smooth form \eqref{eq:passare} then converges unconditionally.
Moreover, also in the complete intersection case, 
a holomorphic function annihilates the Coleff-Herrera product if and only if it is in the ideal $\langle f_1,\ldots,f_q\rangle$;
this last property is called the \emph{duality property} and it was proved independently by Dickenstein-Sessa, \cite{DS},
and Passare, \cite{PDr}.

% By Stokes' theorem the action of $\dbar (1/f)$ on a test form $\varphi$ is given by the limit of
% \begin{equation}\label{eq:HL}
% \lim_{\epsilon\to 0^+} \int_{|f|^2=\epsilon}\varphi/f
% \end{equation}
% as $\epsilon\to 0^+$. It can also be computed using analytic continuation: The Mellin transform, with respect to
% $\epsilon$, of \eqref{eq:HL} equals
% \begin{equation*}
% \int \dbar |f|^{2\lambda}\wedge \varphi/f
% \end{equation*} 
% for $\mathfrak{Re}\, \lambda \gg 0$

\smallskip

In this paper we consider another approach to Coleff-Herrera type products; it is based on analytic continuation
and has been studied in, e.g., \cite{At,Bar,BernGelf,PTcanada,Y}.
For $\lambda_j\in \C$ with $\mathfrak{Re}\, \lambda_j \gg 0$, let
\begin{equation*}
\Gamma_f^{\varphi}(\lambda_1,\ldots,\lambda_q) = 
\int \frac{\dbar |f_q|^{2\lambda_q}\wedge \cdots \wedge \dbar |f_1|^{2\lambda_1}}{f_1\cdots f_q}
\wedge \varphi,
\end{equation*}
where $\varphi$ is a test form.
It is standard to see that $\lambda_1 \mapsto \Gamma_f^{\varphi}(\lambda_1,\ldots,\lambda_q)$ 
has an analytic continuation to a neighborhood of $0$ and that $\Gamma_f^{\varphi}(0,\lambda_2,\ldots,\lambda_q)$ equals
\begin{equation*}
\frac{\dbar |f_q|^{2\lambda_q}}{f_q}\wedge \cdots \wedge \frac{\dbar |f_2|^{2\lambda_2}}{f_2}
\wedge \dbar\frac{1}{f_1} .\, \varphi.
\end{equation*}
From \cite[Proposition~2.1]{AW2} it follows that $\lambda_2 \mapsto \Gamma_f^{\varphi}(0,\lambda_2,\ldots,\lambda_q)$ 
is analytic at $0$, that $\lambda_3 \mapsto \Gamma_f^{\varphi}(0,0,\lambda_3,\ldots,\lambda_q)$ is too, and so on.
% From \cite{AW2} it follows that $\lambda_1 \mapsto \Gamma_f^{\varphi}(\lambda_1,\ldots,\lambda_q)$ is analytic at $0$,
% that $\lambda_2\mapsto \Gamma_f^{\varphi}(0,\lambda_2,\ldots,\lambda_q)$ is analytic at $0$ too, and so on;
% here $\Gamma_f^{\varphi}(0,\lambda_2,\ldots,\lambda_q)$ means the action of $\dbar (1/f_1)$ on 
% $\pm (\dbar |f_q|^{2\lambda_q}/f_q)\wedge \cdots \wedge (\dbar |f_2|^{2\lambda_2})$.
Once one knows that the Coleff-Herrera product is obtained by letting $\epsilon_j\to 0$ successively 
in \eqref{eq:passare} it is not that hard to see that
\begin{equation*}
\dbar \frac{1}{f_q}\wedge \cdots \wedge \dbar \frac{1}{f_1} .\, \varphi=
\Gamma_f^{\varphi}(\lambda_1,\ldots,\lambda_q)|_{\lambda_1=0}\cdots |_{\lambda_q=0},
\end{equation*}
where the expression on the right hand side means that we first let $\lambda_1\to 0$, then let $\lambda_2\to 0$ etc;
see, e.g., \cite[Theorem~2]{LS}.
However, from an algebraic point of view, cf.\ \cite[Theorem 3.2]{BY},
it is often desirable  
to have a current given as the value at $0$ of a single one-variable analytic function; 
this is the motivation for this paper.
From Theorem~\ref{main} below it follows that if $\mu_1 > \cdots > \mu_q >0$ are integers, then 
$\lambda \mapsto \Gamma_f^{\varphi}(\lambda^{\mu_1},\ldots,\lambda^{\mu_q})$,
a priori defined for $\mathfrak{Re}\, \lambda \gg 0$, has an analytic continuation to a neighborhood of
$[0,\infty)\subset \C$ and that the value at $\lambda=0$ equals the Coleff-Herrera product \eqref{eq:CHprod}.
Notice that this way of letting $(\lambda_1,\ldots,\lambda_q)\to 0$ is analogous to limits along admissible paths
in the sense that $\lambda_j$ goes to zero much faster than $\lambda_{j+1}$, $j=1,\ldots,q-1$.

We remark that if $f$ defines a complete intersection then it is showed in \cite{hasamArkiv} that 
$\Gamma_f^{\varphi}(\lambda)$ is analytic in a neighborhood of the half-space $\{\mathfrak{Re}\, \lambda_j\geq 0, j=1,\ldots,q\}$.

\begin{center}
---
\end{center}

Let us now consider a more general setting.
Let $f$ be a section of a Hermitian vector bundle $E$ of rank $m$ over a reduced complex space $X$ of pure dimension $n$. 
In \cite{PTY} and \cite{ABullSci} were introduced currents $U$ and $R$, generalizing the currents $1/f$ and 
$\dbar (1/f)$, respectively. These currents 
are based on Bochner-Martinelli type expressions. To be precise, let
$f=f_{1}e_{1}+\cdots +f_{m}e_{m}$, where $\{e_{k}\}_k$ is a local holomorphic frame for $E$
with dual frame $\{e_k^*\}_k$,
and let $s=s_{1}e_{1}^*+\cdots +s_{m}e_{m}^*$ be the section of the dual bundle $E^*$ with pointwise minimal
norm such that $f\cdot s=|f|^2_{E}$.
For $\lambda\in \C$, $\mathfrak{Re}\, \lambda \gg 0$, we let
\begin{equation}\label{eq:U}
U^{\lambda}:=\sum_{k=1}^{m} |f|^{2\lambda}_{E} \frac{s\wedge (\dbar s)^{k-1}}{|f|^{2k}_{E}},
\end{equation}
where $(0,1)$-forms anticommute with the $e^*_{k}$. It turns out, \cite{ABullSci}, \cite{PTY}, that $\lambda \mapsto U^{\lambda}$,
considered as a current-valued map, has an analytic continuation to a neighborhood of $0$. The value at $\lambda=0$ is a current
$U$ on $X$ that takes values in $\Lambda E^*$; $U$ is the standard extension of $\sum_k s\wedge (\dbar s)^{k-1}/|f|_E^{2k}$
across $\{f=0\}$.
If $E$ has rank $1$, then $U=(1/f)e^*$ for any choice of metric. Let 
\begin{equation}\label{eq:R}
R^{\lambda}:=
1-|f|^{2\lambda}_{E} + \sum_{k=1}^{m} \dbar |f|^{2\lambda}_{E} \wedge 
\frac{s\wedge (\dbar s)^{k-1}}{|f|^{2k}_{E}}.
\end{equation}
Letting $\nabla_{f}:=\delta_{f}-\dbar$, where $\delta_{f}$ denotes interior multiplication with $f$, 
one can check that $R^{\lambda}=1-\nabla_{f}U^{\lambda}$, see \cite{ABullSci} for details.
It follows that $\lambda \mapsto R^{\lambda}$ has an analytic continuation to a neighborhood of $0$ and 
the value at $\lambda=0$ is the current $R$; it is straightforward to check that $R$ has support on $\{f=0\}$. 
If $E$ has rank $1$ then 
$R=\dbar (1/f)\wedge e^*$ and more generally, if $f$ defines a complete intersection then 
$R=\dbar(1/f_{m})\wedge \cdots \wedge\dbar(1/f_{1})\wedge e^*_{1}\wedge \cdots \wedge e^*_{m}$ for any choice of metric,
see \cite{ABullSci} and \cite{PTY}.

The value at $\lambda=0$ of the term $1-|f|^{2\lambda}_{E}$ of $R^{\lambda}$ is the restriction
$\mathbf{1}_{\{f=0\}}$ to the zero set of $f$, see \cite{AW2}. In itself it is zero unless $f$ vanishes identically on some 
components of $X$ in which case it simply is $1$ there. However,
when forming products of $R$'s the role of $\mathbf{1}_{f=0}$ is much more significant, 
cf.\ \cite{ASWY1} and Example~\ref{ex:ASWY}.

Let $f_j$ be a section of a Hermitian vector bundle $E_j$ of rank $m_j$, let $U^j$ and $R^j$ be the 
associated currents, and let $U^{j,\lambda}$ and $R^{j,\lambda}$ be the corresponding $\lambda$-regularizations. 
Following, e.g., \cite{ASWY1} and \cite{LS}
we define products of the $R^j$ recursively as follows. Having defined $R^k\wedge \cdots \wedge R^1$, consider
the current-valued function
\begin{equation*}
\lambda \mapsto
R^{{k+1},\lambda}\wedge R^k\wedge \cdots \wedge R^1,
\end{equation*}
a priori defined for $\mathfrak{Re}\, \lambda \gg 0$. It turns out, see, e.g., \cite{AW2} or \cite{LS}, that it 
can be analytically continued to a neighborhood of $0$, and we let $R^{k+1}\wedge \cdots \wedge R^1$ be the value at $\lambda=0$.

\begin{theorem}\label{main}
Let $\mu_1 > \cdots > \mu_q$ be positive integers. Then the current-valued function
\begin{equation*}
\lambda \mapsto R^{q,\lambda^{\mu_q}}\wedge \cdots \wedge R^{1,\lambda^{\mu_1}},
\end{equation*}
a priori defined for $\mathfrak{Re}\, \lambda \gg 0$, has an analytic continuation to a neighborhood of
the half-axis $[0,\infty)\subset \C$
and the value at $\lambda=0$ is $R^q\wedge \cdots \wedge R^1$.
\end{theorem}

To connect with Coleff-Herrera type products, let $\chi$ be the characteristic function $\mathbf{1}_{[1,\infty)}$ or
a smooth regularization thereof and let 
\begin{equation*}
R^{j,\epsilon_j}:=
1-\chi(|f_j|^{2}_{E_j}/\epsilon_j) + \sum_{k=1}^{m_j} \dbar \chi(|f_j|^{2}_{E_j}/\epsilon_j) 
\wedge \frac{s_j\wedge (\dbar s_j)^{k-1}}{|f_j|^{2k}_{E_j}}.
\end{equation*}
If $\varphi$ is a test form on $X$, then the limit of
\begin{equation}\label{uncond}
\int_X R^{q,\epsilon_q}\wedge \cdots \wedge R^{1,\epsilon_1}\wedge \varphi
\end{equation}
as $\epsilon\to 0$ along an admissible path exists and equals the action of
$R^q\wedge \cdots \wedge R^1$ on $\varphi$, see \cite{LS}. 

\bigskip

Let us mention a version of Theorem~\ref{main} with connection to intersection theory.
Let $f$ be a section of $E$ and let 
\begin{equation*}
M^{\lambda}:= 1-|f|^{2\lambda}_{E} +
\sum_{k\geq 1} \dbar |f|^{2\lambda}_{E} \wedge \frac{\partial \log |f|^{2}_{E}}{2\pi i} \wedge
(dd^c\log |f|^{2}_{E})^{k-1},
\end{equation*}
where $dd^c=\dbar \partial/2\pi i$. It is showed in \cite{ASWY1} that $\lambda \mapsto M^{\lambda}$ has an analytic 
continuation to a neighborhood of $0$ and that the value at
$\lambda=0$ is a positive closed current, which we denote by $M$.
One can give a meaning to the product $(dd^c\log |f|_E^2)^k$ for arbitrary $k$ that extends the classical 
one for $k\leq \textrm{codim}\, \{f=0\}$, and from \cite{ASWY1} it follows that
\begin{equation*}
M = \mathbf{1}_{Z} + \sum_{k\geq 1} \mathbf{1}_{Z} (dd^c\log |f|^{2}_{E})^k,
\end{equation*}
where $\mathbf{1}_{Z}$ is the restriction to the zero set $Z$ of $f$. 
The current $M$ is closely connected to
$R$. For instance, if $X$ is smooth and $D$ is the Chern connection on $E$ then it follows from \cite{AArkiv}
that
\begin{equation*}
M_k=R_k \cdot (Df/2\pi i)^k/k!,
\end{equation*}
where the subscript $k$ means the component of bidegree $(*,k)$.

Let $f_1,\ldots,f_q$ be sections of Hermitian vector bundles $E_j$ and let $M^1,\ldots,M^q$ be the associated current.
One can define products of the $M^{j}$ recursively as for the $R^{j}$ and we have the following
analogue of Theorem~\ref{main}.

\begin{theorem}\label{main2}
Let $\mu_1 > \cdots > \mu_q$ be positive integers. Then the current-valued function
\begin{equation*}
\lambda \mapsto M^{q,\lambda^{\mu_q}}\wedge \cdots \wedge M^{1,\lambda^{\mu_1}},
\end{equation*}
a priori defined for $\mathfrak{Re}\, \lambda \gg 0$, has an analytic continuation to a neighborhood of
the half-axis $[0,\infty)\subset \C$
and the value at $\lambda=0$ is $M^q\wedge \cdots \wedge M^1$.
\end{theorem}

\begin{example}[Example 5.6 in \cite{ASWY1}]\label{ex:ASWY}
Let $\mathcal{J}_x\subset \hol_{X,x}$ be an ideal and let $h_1,\ldots,h_{n}\in \mathcal{J}_x$ be a generic
Vogel sequence of $\mathcal{J}_x$; see, e.g., \cite{ASWY1} for the definition. By the St\"{u}ckrad-Vogel procedure,
\cite{ST}, adapted to the local situation, \cite{Mass1}, \cite{T}, one gets an associated Vogel cycle $V^h$; 
the multiplicities of the components of various dimensions of 
$V^h$ are the Segre numbers, \cite{GG}, used in excess intersection theory. 
By Theorem~\ref{main2} we have that 
\begin{equation*}
\lambda \mapsto  \bigwedge_{k=1}^n
\big( 1-|h_k|^{2\lambda^{\mu_k}}+\dbar|h_k|^{2\lambda^{\mu_k}}\w
\partial \log|h_k|^2/2\pi i\big)
\end{equation*}
is analytic at $0$ and by \cite{ASWY1} the value there is the Lelong current associated with $V^h$;
see \cite{ASWY1} for more details.
\end{example}

\begin{remark}
Assume that $\textrm{codim}\, \cap_j\{f_j=0\} = m_1+\cdots+m_q$.
Then $M^j=(dd^c \log |f_j|^2_{E_j})^{m_j}=[f_j=0]$, where $[f_j=0]$ is the Lelong current of 
the fundamental cycle of $f_j$, and more generally,
\begin{equation*}
M^q\wedge \cdots \wedge M^1=[f_q=0]\wedge \cdots \wedge [f_1=0],
\end{equation*} 
i.e., the current representing the proper intersection of the cycles $[f_j=0]$.

In this case the current-valued function 
\begin{equation*}
(\lambda_1,\ldots,\lambda_q) \mapsto R^{q,\lambda_q}\wedge \cdots \wedge R^{1,\lambda_1}
\end{equation*}
has an analytic continuation to a neighborhood of the origin in $\C^q$, \cite{LS},
and the value at $\lambda=0$ is the $R$-current associated to $\oplus_j f_j$, \cite{W}.
Moreover, by \cite{LS}, \eqref{uncond} depends H\"{o}lder continuously on $\epsilon\in [0,\infty)^q$ if  
$\chi$ is smooth. The smoothness of $\chi$ is necessary in view of the example in
\cite[Section~1]{PTmotex}. 
\end{remark}

\section{Proof of Theorems~\ref{main} and \ref{main2}}
We will actually prove a slightly more general result than Theorem~\ref{main}; we will allow mixed products
of $U^{j}$ and $R^{k}$. Let $P^j$ denote either $U^{j}$ or $R^{j}$ and let $P^{j,\lambda_j}$
be the corresponding $\lambda$-regularization, \eqref{eq:U} or \eqref{eq:R}.
One defines products of the $P^j$ recursively as above.

\medskip

\noindent {\bf Theorem \ref{main}'.}
{\em Let} $\mu_1 > \cdots > \mu_q$ {\em be positive integers. Then the current-valued function}
\begin{equation*}
\lambda \mapsto P^{q,\lambda^{\mu_q}}\wedge \cdots \wedge P^{1,\lambda^{\mu_1}},
\end{equation*}
{\em a priori defined for} $\mathfrak{Re}\, \lambda \gg 0$, 
{\em has an analytic continuation to a neighborhood of the half-axis} $[0,\infty)\subset \C$
{\em and the value at} $0$ {\em is} $P^q\wedge \cdots \wedge P^1$.

\medskip

Let $\pi \colon X'\to X$ be a smooth modification of $X$ such that 
$\{\pi^*f_j=0\}$, $j=1,\ldots,q$, and $\cup_j\{\pi^*f_j=0\}$ are normal crossings divisors. 
Then locally in $X'$ we can write $\pi^*f_j=f_j^0 f_j'$, where $f_j^0$ is a monomial in local coordinates
and $f_j'$ is a non-vanishing holomorphic tuple. It follows that $s_j = \bar{f}_j^0 s'_j$, where
$s'_j$ is a smooth section. A straightforward computation shows that
\begin{equation*}
\pi^* R^{j,\lambda_j}= 1-|f_j^0|^{2\lambda_j} u_j^{2\lambda_j} +
\sum_{k=1}^{m_j} \frac{\dbar(|f_j^0|^{2\lambda_j} u_j^{2\lambda_j})}{(f_j^0)^k}\wedge \vartheta_{jk},
\end{equation*}
\begin{equation*}
\pi^* U^{j,\lambda_j}= 
\sum_{k=1}^{m_j} \frac{|f_j^0|^{2\lambda_j} u_j^{2\lambda_j}}{(f_j^0)^k}\wedge \vartheta_{jk},
\end{equation*}
where $u_j$ is a smooth non-vanishing function $\vartheta_{jk}$ is a smooth form.
In the same way, %Moreover, from \cite[Section~4]{ASWY1}
\begin{equation*}
\pi^* M^{f_j,\lambda_j}= 1-|f_j^0|^{2\lambda_j} u_j^{2\lambda_j} +
\sum_{k\geq 1} \dbar\big(|f_j^0|^{2\lambda_j} u_j^{2\lambda_j}\big) \wedge \partial \log(|f_j^0|^2 u_j^2)\wedge \omega_{jk},
\end{equation*}
where $\omega_{jk}$ is smooth, cf.\ \cite[Section~4]{ASWY1}.
Theorems \ref{main}' and \ref{main2} are immediate consequences of the following 
quite technical lemma; indeed $\partial \log(|f_j^0|^2 u_j^2) = df_j^0/f_j^0 + 2\partial u_j/u_j$.

\begin{lemma}\label{pilligt}
Let $u_1,\ldots,u_r$ be smooth non-vanishing functions defined in some neighborhood of the origin in $\C^n$, with coordinates 
$x_1,\ldots, x_n$. For $\lambda=(\lambda_1,\ldots, \lambda_r)\in\C^r$, $\mathfrak{Re}\,\lambda_j \gg 0$, 
$\alpha_1,\ldots, \alpha_r\in\mathbb{N}^n$, and $k_1,\ldots,k_r\in \mathbb{N}$, let  
\begin{equation*}
\Gamma (\lambda):= \frac{|u_rx^{\alpha_r}|^{2\lambda_r}\cdots |u_{p+1}x^{\alpha_{p+1}}|^{2\lambda_{p+1}}
 \bar{\partial}|u_px^{\alpha_p}|^{2\lambda_p}\wedge \cdots \wedge
\bar{\partial}|u_1x^{\alpha_1}|^{2\lambda_1}}{x^{k_r\alpha_r}\cdots x^{k_1\alpha_1}};
\end{equation*}
here $x^{k_{\ell}\alpha_\ell}=x_1^{k_{\ell}\alpha_{\ell,1}}\cdots x_n^{k_{\ell}\alpha_{\ell,n}}$ if 
$\alpha_\ell=(\alpha_{\ell,1},\ldots, \alpha_{\ell,n})$. 
If $\sigma$ is a permutation of $\{1,\ldots, r\}$, write  
$
\Gamma^{\sigma}(\lambda_1,\ldots,\lambda_r):=\Gamma(\lambda_{\sigma(1)},\ldots,\lambda_{\sigma(r)})\,.
$ 
%Then $\Gamma^{\sigma}$ has an analytic continuation to $\Re\lambda_1>-\epsilon_1$ for some $\epsilon_1$. The value at $\lambda=0$, which we denote by $\Gamma^{\sigma}|_{\lambda_1=0}$, has an analytic continuation to $\Re\lambda_2> -\epsilon_2$, etc. 

Let $\mu_1,\ldots, \mu_r$ be positive integers. % and let $\kappa\in\C$, $\Re\kappa >> 0$. 
Then $\Gamma^{\sigma}(\kappa^{\mu_1},\ldots,\kappa^{\mu_r})$  has an analytic continuation to a connected 
neighborhood of the half-axis $[0,\infty)$ in $\C$,
and if  $\mu_1>\ldots>\mu_r$, then
\begin{equation}\label{hovas}
\Gamma^{\sigma}(\kappa^{\mu_1},\ldots,
\kappa^{\mu_r})\mid_{\kappa=0}=
\Gamma^{\sigma}(\lambda_1,\ldots, \lambda_r)\mid_{\lambda_1=0}\cdots
\mid_{\lambda_r=0}.
\end{equation} 
%%where the right hand side is well-defined. 
%where the right hand side of 
%(in the right-hand side of (\ref{hovas}), we mean that we take the analytic continuation in $\lambda$ 
%from a truncated angular sector $S_{\epsilon,R(\epsilon)}$ (cf., (\ref{sector})) and follow it up to $\lambda=0$). 
\end{lemma}

The reason for the permutation $\sigma$ is that we have mixed products of $U$'s and $R$'s in Theorem~\ref{main}'.

\begin{proof}
To begin with let us assume that all $u_j=1$.
A straightforward computation shows that
\begin{equation*}
\Gamma(\lambda) =
\lambda_1 \cdots \lambda_p
\frac{\prod_{j=1}^r|x^{\alpha_j}|^{2\lambda_j}}{x^{\sum_{j=1}^r k_j\alpha_j}}
\sum'_{I} 
A_I
\frac{d\bar{x}_{i_1}\wedge \cdots \wedge
d\bar{x}_{i_p}}{\bar{x}_{i_1}\cdots \bar{x}_{i_p}} =:\lambda_1\cdots\lambda_p\sum'_I \Gamma_I,
\end{equation*}
where the sum is over all increasing multi-indices $I=\{i_1,\ldots, i_p\} \subset \{1,\ldots,n\}$ and 
$A_I$ is the determinant of the matrix $(\alpha_{\ell,i_j})_{1\leq\ell\leq p, 1\leq j\leq p}$.  

Pick a non-vanishing summand $\Gamma_I$; without loss of generality, assume that $I=\{1,\ldots, p\}$ and $A_I=1$. 
With the notation $b_k(\lambda):=\sum_{\ell=1}^r\lambda_\ell \alpha_{\ell,k}$ for $1\leq k\leq n$,
\begin{multline*}
\Gamma_I=\frac{\prod_{k=1}^n|x_k|^{2b_k(\lambda)}}{x^{\sum_{j=1}^r k_j\alpha_j}}
\frac{d\bar{x}_{1}\wedge \cdots \wedge d\bar{x}_{p}}{\bar{x}_{1}\cdots
\bar{x}_{p}} 
= \\
\frac{1}{b_1(\lambda)\cdots b_p(\lambda)}
\frac{\bigwedge_{k=1}^p\bar{\partial}|x_k|^{2b_k(\lambda)}\prod_{k=p+1}^n|x_k|^{2b_k(\lambda)}}{x^{\sum_{j=1}^r
k_j\alpha_j}}.
\end{multline*} 
%It is well known (and easy to check) that the current-valued function
Now the current-valued function 
$$
\widetilde{\Gamma}_I:(\lambda_1,\ldots,\lambda_r)
\mapsto
\frac{\bigwedge_{j=1}^p\bar{\partial}|x_j|^{2b_j(\lambda)}\prod_{j=p+1}^n|x_j|^{2b_j(\lambda)}}{x^{\sum k_j\alpha_j}}
$$
has an analytic continuation to a neighborhood of the origin in
$\mathbb{C}^r$; in fact, it is a tensor product of one-variable currents. In particular,  
$
\widetilde{\Gamma}_I(\kappa^{\mu_1},\ldots, \kappa^{\mu_r})\mid _{\kappa=0}=
\widetilde{\Gamma}_I(\lambda)\mid_{\lambda_1=0}\cdots\mid_{\lambda_r=0}
$.
Let 
$$
\gamma (\lambda)=\frac{\lambda_1\cdots\lambda_p}{b_1(\lambda)\cdots b_p(\lambda)}
$$ 
and 
$\gamma^\sigma=\gamma(\lambda_{\sigma(1)},\ldots, \lambda_{\sigma(r)})$. 
We claim that if $\mu_1>\ldots >\mu_r$, then 
$$
\gamma^\sigma(\lambda)\mid_{\lambda_1=0}\cdots \mid_{\lambda_r=0}=
\gamma^{\sigma}(\kappa^{\mu_1},\ldots,\kappa^{\mu_r})|_{\kappa=0},
$$ 
where it is a part of the claim that both sides make sense.
%Thus \eqref{product} follows in this case. 

Let us prove the claim. Since $A_I =1$, reordering  the factors $b_1,\ldots, b_p$ and 
multiplying $\gamma(\lambda)$ by a non-zero constant, we may assume that $\alpha_{kk}=1$, $k=1,\ldots, p$, so that 
\begin{equation*}
\gamma(\lambda)=
\frac{\lambda_1}{\lambda_1+\alpha_{21}\lambda_2+\cdots +\alpha_{r1}\lambda_r}
\cdots
\frac{\lambda_p}{\alpha_{p1}\lambda_1+\cdots + \lambda_p+\cdots +\alpha_{rp}\lambda_r}.
\end{equation*}
For $j<r$ set $\tau_j:=\lambda_j/\lambda_{j+1}$ and $\widetilde\gamma^\sigma(\tau_1,\ldots, \tau_{r-1}):=\gamma^\sigma(\lambda)$; 
notice that $\gamma^\sigma$ is $0$-homogeneous, so that $\widetilde\gamma^\sigma$ is well-defined. 
%Moreover, note that $\lambda_j=\tau_j\cdots\tau_{r-1}\lambda_r$. It follows that $\widetilde\gamma^\sigma$ 
%%consists of $p$ factors of the form 
Then $\lambda_j=\tau_j\cdots\tau_{r-1}\lambda_r$, and therefore   $\widetilde\gamma^\sigma$ consists of $p$ factors of the form
\begin{equation}\label{leonard}
%\frac{\lambda_k}
%{b_1\lambda_1+\cdots+\lambda_k+\cdots+b_r\lambda_r}
%=\\
\frac{\tau_k\cdots\tau_{r-1}}
{\alpha_{k1}\tau_1\cdots\tau_{r-1}+\cdots+\tau_k\cdots\tau_{r-1}+\cdots+\alpha_{k,r-1}\tau_{r-1}+\alpha_{kr}}.
\end{equation}
Observe that \eqref{leonard} is holomorphic in $\tau$ in some neighborhood of the origin. %%% say $|\tau_j|< \epsilon$. 
Indeed, if $\alpha_{kr}\neq 0$, then \eqref{leonard} is clearly holomorphic, whereas if $\alpha_{kr}=0$ we can factor 
out $\tau_{r-1}$ from the denominator and numerator. In the latter case %if $b_{r-1}\neq 0$
 \eqref{leonard} is clearly holomorphic if $\alpha_{k,r-1}\neq 0$ etc; since $\alpha_{kk}=1$ this procedure 
eventually stops. 
%By repating we see that \eqref{leonard} is holomorphic in $\tau$ in some 
%%neighborhood of the origin, say $|\tau_j|\leq \epsilon$. If $b_r\neq 0$,
%%\eqref{leonard} is clearly holomorphic when $|\tau_1|,\ldots,|\tau_{r-1}|$ 
%%are small. If $b_r=0$ we can factor our $\tau_{r-1}$ from the denominator and numerator. 
%%By repating we see that \eqref{leonard} is holomorphic in $\tau$ in some neighborhood of the origin, say $|\tau_j|\leq \epsilon$. 
%
Hence, $\tilde{\gamma}^{\sigma}(\tau)$ is holomorphic in a neighborhood of $0$.
It follows that $\gamma^\sigma(\kappa^{\mu_1},\ldots, \kappa^{\mu_r})=
\widetilde\gamma^\sigma(\kappa^{\mu_1-\mu_2},\ldots,\kappa^{\mu_{r-1}-\mu_r})$ is holomorphic in a neighborhood 
of $0$ and since the denominator of $\gamma^\sigma(\kappa^{\mu_1},\ldots, \kappa^{\mu_r})$ is a polynomial 
in $\kappa$ with non-negative coefficients it is in fact holomorphic in a neighborhood  
of $[0,\infty)$. Moreover, $\gamma^\sigma(\lambda_1,\ldots, \lambda_r)$ is holomorphic in 
$\Delta=\{|\lambda_1/\lambda_2|< \epsilon,\ldots, |\lambda_{r-1}/\lambda_r|< \epsilon\}$. 
Let us now fix  $\lambda_2\neq 0,\ldots, \lambda_r\neq 0$ in $\Delta$. Then $\gamma^\sigma(\lambda)$ is 
holomorphic in $\lambda_1$ in a neighborhood of the origin. Next, 
for $\lambda_3\neq 0,\ldots, \lambda_r\neq 0$ fixed in $\Delta$, $\gamma^\sigma(\lambda)|_{\lambda_1=0}$ 
is holomorphic in $\lambda_2$ in a neighborhood of the origin, etc. It follows that 
\begin{equation*}
\gamma^\sigma(\lambda)|_{\lambda_1}\cdots |_{\lambda_r=0}=\widetilde\gamma^\sigma(\tau)|_{\tau=0}=
\gamma^\sigma(\kappa^{\mu_1},\ldots, \kappa^{\mu_r})|_{\kappa=0},
\end{equation*}
which proves the claim. Thus \eqref{hovas} follows in the case $u_j=1$, $j=1,\ldots,r$.

%\begin{equation*}
%A\lambda_1\cdots \lambda_p \frac{|x_1^{(\lambda,a_1)}|^2\cdots
%|x_n^{(\lambda,a_n)}|^2}{x^{\sum \alpha_j}}
%\frac{d\bar{x}_{1}\wedge \cdots \wedge d\bar{x}_{p}}{\bar{x}_{1}\cdots
%\bar{x}_{p}} \hspace{3cm}
%\end{equation*}
%\begin{equation*}
%=A\frac{\lambda_1\cdots \lambda_p}{(\lambda,a_1)\cdots (\lambda,a_p)}
%\frac{\bigwedge_1^p\bar{\partial}|x_j|^{2(\lambda,a_j)}\prod_{p+1}^n|x_j|^{2(\lambda,a_j)}}{x^{\sum
%\alpha_j}}.
%\end{equation*}
%It is well known (and easy to check) that the current-valued function
%$(\lambda_1,\ldots,\lambda_r)
%\mapsto
%\bigwedge_1^p\bar{\partial}|x_j|^{2(\lambda,a_j)}\prod_{p+1}^n|x_j|^{2(\lambda,a_j)}/x^{\sum
%\alpha_j}$
%has an analytic continuation to a neighborhood of the origin in
%$\mathbb{C}^r$. Thus, the lemma follows
%from the elementary fact that, if
%\begin{equation*}
%\gamma (\lambda_1,\ldots,\lambda_r)=\frac{\lambda_1\cdots
%\lambda_p}{(\lambda,a_1)\cdots (\lambda,a_p)}
%\end{equation*}
%and $\det (a_{ij})_{1\leq i,j \leq p}\neq 0$, then
%\begin{equation*}
%\gamma^{\sigma}\mid_{\lambda_1=0}\cdots \mid_{\lambda_r=0}=\lim_{\lambda
%\to 0}
%\gamma^{\sigma}(\lambda^{\mu_1},\ldots,\lambda^{\mu_r}).
%\end{equation*}

\smallskip

Now, consider the general case. 
Replace each $|u_j|^{2\lambda_j}$ in $\Gamma(\lambda)$ by $|u_j|^{2\omega_j}$, where $\omega_j\in \C$.
Then $\Gamma$ is a sum of terms of the following representative form:
\begin{equation}\label{lurk}
\prod_{j=p+1}^r |u_j|^{2\omega_j} \prod_{j=1}^{p'} |u_j|^{2\omega_j} \bigwedge_{p'+1}^p \dbar |u_j|^{2\omega_j} \wedge
\frac{\prod_{j=p'+1}^r |x^{\alpha_j}|^{2\lambda_j} \bigwedge_{j=1}^{p'} \dbar |x^{\alpha_j}|^{2\lambda_j}}{
x^{k_r\alpha_r}\cdots x^{k_1 \alpha_1}}
\end{equation}
Fixing all $\lambda_j$ and $\omega_j$ except for $\lambda_{\sigma(1)}$ and $\omega_{\sigma(1)}$, \eqref{lurk}
becomes an analytic (current-valued) function $g(\lambda_{\sigma(1)},\omega_{\sigma(1)})$ in a neighborhood of $0\in \C^2$. 
Thus, the value at $0$ of $g(\lambda_{\sigma(1)},\lambda_{\sigma(1)})$ is the same as first letting $\omega_{\sigma(1)}=0$ 
(which corresponds to setting $u_{\sigma(1)}=1$) and then
letting $\lambda_{\sigma(1)}=0$ in $g(\lambda_{\sigma(1)},\omega_{\sigma(1)})$.
Continuing analogously for $(\lambda_{\sigma(2)},\omega_{\sigma(2)})$ and so on, it follows that the right hand side of 
\eqref{hovas} is independent of the $u_j$. 

To see that the left hand side of \eqref{hovas} 
is independent of $u_j$, replace each $\lambda_j$ in \eqref{lurk} by $\kappa^{\mu_{\sigma(j)}}$
and denote the resulting expression by $\tilde{g}(\kappa,\omega_1,\ldots,\omega_r)$. Then $\tilde{g}$ is clearly analytic 
in the $\omega_j$ and by the first part of the proof it is also analytic in a neighborhood of $[0,\infty)\subset \C_{\kappa}$.
Hence, $\tilde{g}$ is analytic in a neighborhood of $0\in \C^{r+1}$. The left hand side of \eqref{hovas}
is obtained by evaluating $\kappa \mapsto \tilde{g}(\kappa,\kappa^{\mu_{\sigma(1)}},\ldots,\kappa^{\mu_{\sigma(r)}})$
at $\kappa =0$; this is thus the same as evaluating $\tilde{g}(\kappa,0)$ (which corresponds to setting all $u_j=1$) at $\kappa=0$.
Hence also the left hand side of \eqref{hovas} is independent of the $u_j$ and the lemma follows.

% $u_j^{2\lambda_j}$ by $u_j^{2\omega_j}$ 
% in $\Gamma(\lambda)$. Then, by arguments as above, the function
% $(\kappa,\omega_1,\ldots, \omega_r) \mapsto 
% \Gamma^\sigma(\kappa^{\mu_1},\ldots,\kappa^{\mu_r})$ 
%\Gamma^\sigma(\kappa^{\mu_1},\ldots,\kappa^{\mu_r},\omega_1,\ldots,\omega_r)$ 
% is holomorphic in a neighborhood of the origin in
% $\C_\lambda \times \C^r_\omega$  since it is analytic in each variable. In particular, $\Gamma(\kappa^{\mu_1},\ldots,
% \kappa^{\mu_r})\mid_{\kappa=0}$ is obtained by first setting each $\omega_j=0$, which corresponds to setting $u_j=1$ and 
% thus brings us back to the special case treated above. 
% This completes the proof of the lemma.  
\end{proof}

\end{document}